\documentclass[review]{elsarticle}

\usepackage{lineno,hyperref}
\usepackage{graphicx}
\usepackage{amssymb}
\usepackage{epstopdf}
\usepackage{epsfig}
\usepackage{amsmath}
\usepackage{amsthm}
\usepackage{caption}

\newtheorem{theorem}{Theorem}
\newtheorem{lemma}{Lemma}
\newtheorem{remark}{Remark}
\newtheorem{example}{Example}
\journal{Chaos, Solitons and Fractals}

\begin{document}
\captionsetup[figure]{labelfont={bf},labelformat={default},labelsep=period,name={Fig.}}
\begin{frontmatter}

\title{A criterion for asymptotic stability of general fractional-order linear time-invariant systems with incommensurate orders}

\author{Jing Yang}
\address{University of Electronic Science and Technology of China, Chengdu 611731,
China}

\author{Xiaorong Hou$^*$}
\address{University of Electronic Science and Technology of China, Chengdu 611731,
China}
\cortext[mycorrespondingauthor]{Corresponding author}
\ead{Houxr@uestc.edu.cn}

%

%

\begin{abstract}
A criterion on the asymptotic stability of fractional-order systems with incommensurate orders is proposed in this paper. Existing methods always assume order parameters be rational numbers or the ratios of any two orders be rational numbers. In engineering applications, order parameters are more likely to be uncertain which may be real numbers. Furthermore, the boundary of the stable parameter region is determined, which decomposes parameter space into the finite number of connected regions. All systems whose parameters belong to the same region have the same stability. Each region only needs checking one point to determine the stability of the region. The method established in this paper involves low computational complexity and clearly gives the relationship between order parameters and stability. Some examples show the advantages of this method.
\end{abstract}

\begin{keyword}
Fractional-order system\sep Incommensurate orders\sep Stability \sep Uncertain system\sep Parameter space
\end{keyword}

\end{frontmatter}


\section{Introduction}
Fractional-order systems have been widely studied because they can model real systems more accurately \cite{liu2004improved,peaucelle2000new}. With the development of fractional calculus theories, researchers have found that many systems from engineering applications are of fractional order. For example, the results of biological neurons show that fractional-order neural network models are more practical \cite{kaslik2011dynamics,kaslik2012nonlinear}. Researchers have realized these facts and discussed some problems by using fractional-order systems. Among these problems, the order parameter is always an important subject which makes systems very different from integer-order systems \cite{mujumdar2015observer}. Some results on fractional-order systems with commensurate orders and incommensurate orders have been proposed in recent years \cite{deng2007stability,tavazoei2008chaotic,chang2010chaos,daftardar2010chaos,diethelm2017asymptotic,brandibur2021stability,huseynov2021explicit,brandibur2021exact,djenina2020stability,brandibur2017stability}.
\par
Stability is a basic problem when we consider analyzing dynamic problems of systems. For a fractional-order linear time-invariant system with order $0<\alpha<2$, the eigenvalues $\lambda$ of the system matrix determine the stability of the system. The condition $\left | arg(\lambda) \right |>\frac{\alpha \pi }{2}$ of all eigenvalues $\lambda$ means that the system is asymptotically stable \cite{matignon1996stability,radwan2009stability}. Most methods on the stability of fractional-order systems with commensurate orders are based on this condition \cite{zhang2021exact,sabatier2010lmi,jiao2012robust,lu2013stability,wei2017completeness,2015Robust,2020Output}. Fractional-order systems with incommensurate orders are more general than systems with commensurate orders. Some scientists discussed that systems in 2007 and gave the stability condition that all eigenvalues $\lambda$ satisfy $\left | arg(\lambda) \right |>\frac{\pi }{2m}$, where incommensurate orders are rational numbers and $m$ is the lowest common multiple(LCM) of denominators of orders \cite{deng2007stability}. Several classical fractional-order systems are of incommensurate orders, such as Bagley-Torvik system \cite{wang2010general} and Basset force system \cite{baleanu2013fractional}, which can be directly analyzed with the stability condition \cite{deng2007stability}. The above stability results have been extended to fractional-order nonlinear systems \cite{tavazoei2008chaotic,chang2010chaos,daftardar2010chaos}.
\par
Fractional-order systems with multiple orders are more practical. For example, since electronic components, such as capacitors and inductors, are of fractional-order, it is more accurate to model power electronic devices as multi-order fractional-order systems \cite{bohannan2006electrical,bohannan2002analog}. The n-dimensional fractional-order system is more likely to be with incommensurate orders. Researchers have modeled some real systems as fractional-order systems and analyzed their dynamic problems in the case of incommensurate orders \cite{deng2007stability,tavazoei2008chaotic,chang2010chaos,daftardar2010chaos}. Researchers have paid more attention to systems whose incommensurate orders are all rational numbers. Let the LCM of denominators of all orders be a new parameter, so that systems with incommensurate orders can be rewritten as systems with commensurate orders. In real systems, it is hard to determine the values of order parameters and whether these orders are rational numbers or not. When we consider these order parameters as uncertain parameters, it means that order parameters may be real numbers including rational numbers and irrational numbers. Most papers fixed incommensurate orders as rational numbers to obtain some results. In recent years, some papers considered some special fractional-order systems with two or three different order parameters \cite{vcermak2015asymptotic,diethelm2017asymptotic,brandibur2021stability,brandibur2021exact}, and proposed some methods to solve stable parameter region and unstable parameter region in case of order-dependent or order-independent, respectively. In these results, two- and three-dimensional systems were discussed. The relationship between order parameters and stability for the n-dimensional fractional-order systems is still a problem that need to be analyzed.
\par
For fractional-order systems with incommensurate orders, the following problems need to be considered: First, if order parameters are rational numbers with large denominators, the LCM of denominators will be so large that there are too many roots of characteristic equations of systems. A large number of roots will have high computational complexity. Second, the existing methods use the numerical calculation method to describe the stable region and unstable region. Although these methods are applicable for some low dimension systems, the boundary of the stable region is not clear and unproven.
\par
In this paper, we consider the problem of stability of the general fractional-order linear time-invariant systems with incommensurate orders, arbitrary real numbers. Based on strict proof, a criterion for asymptotic stability of fractional-order systems with incommensurate orders is established, where the order parameters are real numbers or uncertain values. Based on the criterion, the boundary of the stable parameter region of the systems with multiple orders is determined in parameter space, which decomposes parameter space into the finite number of connected regions. All systems whose parameters belong to the same region have the same stability. Each region only needs checking one point to determine the stability of the region. The method proposed in this paper analyzes fractional-order systems with incommensurate orders directly and considers the stability problem in parameter space.
\section{Problem Formulation and Preliminaries}
Consider the following fractional-order system with incommensurate orders:
\begin{equation}
\textrm{D}^{\alpha}x\left ( t \right )=Ax\left ( t \right ),
\end{equation}
where $\alpha=\left [ \alpha_{1},\alpha_{2},\cdots ,\alpha_{n} \right ]^{T}$ and all $\alpha_{i}(i=1,2,\cdots,n)$ are real numbers between 0 and 2. Since the stability of linear systems is independent of initial conditions, initial conditions are omitted in this paper.
\par Its representation in components is:
\begin{equation}
\left\{\begin{matrix}
\textrm{D}^{\alpha_{1}}x_{1}\left ( t \right )=a_{11}x_{1}\left ( t \right )+a_{12}x_{2}\left ( t \right )+\cdots+a_{1n}x_{n}\left ( t \right )\\
\textrm{D}^{\alpha_{2}}x_{2}\left ( t \right )=a_{21}x_{2}\left ( t \right )+a_{22}x_{2}\left ( t \right )+\cdots+a_{2n}x_{n}\left ( t \right )\\
\cdots \cdots \\
\textrm{D}^{\alpha_{n}}x_{n}\left ( t \right )=a_{n1}x_{1}\left ( t \right )+a_{n2}x_{2}\left ( t \right )+\cdots+a_{nn}x_{n}\left ( t \right )
\end{matrix}\right..
\end{equation}
\par The characteristic equation of system (1) is $\Delta _{(\alpha ,A)}\left ( s \right )=0$, where
\begin{equation}\label{eq3}
\Delta _{(\alpha ,A)}\left ( s \right )=det\left ( diag\left ( s^{\alpha_{1}},s^{\alpha_{2}},\cdots ,s^{\alpha_{n}} \right )-A \right ).
\end{equation}
\begin{lemma} \cite{petravs2011fractional}
For system (1), if all $\alpha_{i}(i=1,2,\cdots,n)$ are rational numbers between 0 and 2, then the system is asymptotically stable if and only if all the roots $\lambda $ of the characteristic equation
\begin{equation}\label{eq4}
det\begin{pmatrix}
\lambda^{m\alpha_{1}}-a_{11}&-a_{12}&\cdots& -a_{1n}\\
-a_{21}&\lambda^{m\alpha_{2}}-a_{22}&\cdots& -a_{2n}\\
\vdots & & \vdots \\
-a_{n1}&-a_{n2}&\cdots &\lambda^{m\alpha_{n}}-a_{nn}
\end{pmatrix}=0,
\end{equation}
satisfy
\begin{equation}\label{eq5}
\left | arg\left ( \lambda \right ) \right |>\frac{\pi }{2m},
\end{equation}
where $m$ is the LCM of the denominators $u_{i}$ of $\alpha_{i}$, $\alpha_{i}=v_i/u_i\in(0,2)$, $v_i,u_i\in \mathbb{Z}^{+}$ for $i=1,2,\cdots ,n$.
\end{lemma}
\par
Based on Lemma 1, we have $mq_1, mq_2, \cdots, mq_n\in \mathbb{Z}^+$. The number of roots of Eq.(4) is $mq_1+mq_2+ \cdots+mq_n$, which means we need to check $mq_1+mq_2+ \cdots+mq_n$ roots whether each root satisfies $\left | arg\left ( \lambda \right ) \right |>\frac{\pi }{2m}$. For example, consider the following system:
\begin{equation}\label{eq6}
\left\{\begin{matrix}
_{0}\textrm{D}_{t}^{q_{1}}x_{1}\left ( t \right )=a_{11}x_{1}\left ( t \right )+a_{12}x_{2}\left ( t \right )\\
_{0}\textrm{D}_{t}^{q_{2}}x_{2}\left ( t \right )=a_{21}x_{2}\left ( t \right )+a_{22}x_{2}\left ( t \right )
\end{matrix}\right..
\end{equation}
Its characteristic equation is:
\begin{equation}\label{eq7}
\lambda^{mq_1+mq_2}-a_{22}\lambda^{mq_1}-a_{11}\lambda^{mq_2}+a_{11}a_{22}-a_{12}a_{21}=0.
\end{equation}
If $q_1=67/100, q_2=81/100$, then $mq_1+mq_2=148$. In engineering applications, the order parameters of these systems are usually complex, the values of order parameters are more likely to be 0.985, 0.9776 and so on. In real applications, it is difficult to analyze systems by using Lemma 1 directly. Especially, If these order parameters are uncertain, then Lemma 1 can not be used to analyze the stability of the system.
\par
Fractional-order systems with interval order are often discussed in recent years \cite{liao2011robust,li2012robust}. Some scientists considered systems with interval order and gave some sufficient conditions on the stability of systems. Interval order includes rational numbers and irrational numbers. The existing methods are hard to obtain complete results.
\par
Based on Lemma 1, if all $\alpha_{i}$ are rational numbers, Eq.(3) can be transformed into an integer-order polynomial equation as Eq.(4). When $\alpha=\alpha_1=\alpha_2=\cdots=\alpha_n$ \cite{deng2007stability}, denoted $\lambda^m$ by $s$ in Eq.(4), we get the Eq.(3), where $\left | arg\left ( s \right ) \right |=\left | arg\left ( \lambda^m \right ) \right|$.
\begin{lemma} \cite{matignon1996stability}
For system (1), if $\alpha_{1}=\alpha_{2}=\cdots =\alpha_{n}=\alpha$, then the characteristic equation becomes $det(s^\alpha I-A)=0$. In this case, system (1) is stable if and only if each eigenvalue $\lambda$ of matrix $A$ satisfies $|arg(\lambda)|>\frac{\alpha \pi }{2}$.
\end{lemma}
\par
Most of the existing researches considered some particular systems with multiple orders, and used Lemma 1 to analyze the dynamic problems. Based on Lemma 1, for system (1) with rational number orders, the critical line in the complex plane is the set of $z \in \mathbb{C}$ satisfying $\left |arg(z)  \right |=\frac{\pi }{2m}$. Based on Lemma 2, for system (1) with commensurate orders, the critical line in the complex plane is the set of $z \in \mathbb{C}$ satisfying $\left |arg(z)  \right |=\frac{\alpha \pi }{2}$. However, for general system (1), The critical line is unknown.
\par
In this paper, we consider the stability of the general n-dimensional system.
\section{Main Results}
In this section, we give a criterion of asymptotic stability on system (1).
\begin{lemma}\cite{kilbas2006theory}
The solution of system (1) continuously depends on parameters $\left ( \alpha ,A \right )$.
\end{lemma}
\proof
Let $\Delta _{(\alpha ,A)}\left ( s \right )=\sum_{\alpha }a_{\alpha }s^{\alpha }$, the system (1) can be rewritten as follows:
\begin{equation}
\sum_{\alpha}a_{\alpha}\textrm{D}_t^{\alpha}x(t)=0.
\end{equation}
Based on the Theorem 5.14 in \cite{kilbas2006theory}, the solution of system(8) continuously depends on parameters $\left ( \alpha ,a_{\alpha} \right )$. System (8) is equivalent to system (1), it is obvious that the solution of system (1) continuously depends on parameters $\left ( \alpha ,A \right )$.
\endproof
\par
Let $x\left ( t;\alpha ,A \right )$ be the solution of system (1) with parameters $\left ( \alpha ,A \right )$. Based on the definition of asymptotic stability, we call the point $\left ( \alpha ,A \right )$ is a stable parameter point if $\exists (\alpha,A)\in \mathbb{R}^{n}\times \mathbb{R}^{n\times n}$ such that $x(t;\alpha ,A)\rightarrow 0$ when $t\rightarrow \infty$. The set of stable parameter points is the stable parameter region of system(1), denoted by $S\left ( \alpha ,A \right )$.
\begin{lemma}
$S\left ( \alpha ,A \right )$ is an open set.
\end{lemma}
\proof
To show that $S\left ( \alpha ,A \right )$ is open, take any point $p\in S\left ( \alpha ,A \right )$ and show that every point in a neighborhood  of $p$ belongs to $S\left ( \alpha ,A \right )$. To the end, let $T>0$ be large enough that $\left \| x\left ( T;p \right ) \right \|<\epsilon /2$.
\par Consider the neighborhood $\left \|q-p \right \|<\delta $ of $p$. By continuous dependence of the solution on parameters $\left ( \alpha ,A \right )$, we can choose $\delta$ small enough to ensure that for any point $q$ in the neighborhood $\left \| q-p \right \|<\delta $, the solution at time $T$ satisfies
\begin{equation}
\left \| x(T;p)-x(T;q) \right \|<\epsilon /2,
\end{equation}
then
\begin{equation}
\left \| x(T;q)\right \|\leq \left \| x(T;q)-x(T;p) \right \|+\left \| x(T;p) \right \|<\epsilon.
\end{equation}
This shows that the point $q$ is inside $S(\alpha,A)$. Thus, the set $S(\alpha,A)$ is open.
\endproof
\begin{theorem}
For system (1), the stable parameter region and its boundary are described as follows:
\begin{equation}
S\left ( \alpha ,A \right )=\left \{ (\alpha ,A) \in \mathbb{R}^{n}\times \mathbb{R}^{n\times n}|\left | arg(\lambda_{(\alpha ,A)}) \right |>\frac{\pi }{2}\right \},
\end{equation}
\begin{equation}
R\left ( \alpha ,A \right )=\left \{ (\alpha ,A) \in \mathbb{R}^{n}\times \mathbb{R}^{n\times n}|\left | arg(\lambda_{(\alpha ,A)}) \right |=\frac{\pi }{2}\right \},
\end{equation}
where $\lambda_{(\alpha ,A)}$ is the root of $\Delta_{(\alpha ,A)}(s)=0$. If $\left | arg(\lambda_{(\alpha ,A)}) \right |=\frac{\pi }{2}$, then the $\lambda_{(\alpha ,A)}$ is called the critical root of system (1).
\endproof
\end{theorem}
\proof
\par 1. Consider that all $\alpha _{i}\in \left ( 0,2 \right )$ are rational numbers for $i=1,2,\cdots ,n$. By Lemma 1, denoting $\lambda^m$ by $\lambda_{(\alpha ,A)}$, we have $\left | arg\left ( \lambda_{(\alpha ,A)} \right ) \right |=\left | arg\left ( \lambda^m \right ) \right |$. Then the conclusion holds.
\par 2. Consider that all $\alpha _{i}\in \left ( 0,2 \right )$ are not all rational numbers for $i=1,2,\cdots ,n$.
\par 2.1 Suppose $\exists (\alpha ,A)\in R{(\alpha ,A)}$, $\left | arg\left ( \lambda_{(\alpha ,A)} \right ) \right |=r>\pi /2$. $S{(\alpha ,A)}$ is open set from Lemma 4, so $(\alpha ,A)\notin S{(\alpha ,A)}$. Thus, $\exists \varepsilon >0$, $\forall T>0$, $\exists t>T$, such that
\begin{equation}
\left \| x(t;\alpha ,A) \right \|>\varepsilon.
\end{equation}
Because $x(t;\alpha ,A)$ continuously depends on parameters $(\alpha,A)$, there exists a vector $\bar{\alpha}$ whose components are rational numbers, such that
\begin{equation}
\pi/2<\left | arg(\lambda_{(\alpha ,A))} \right |<r,(\alpha ,A)\in S{(\alpha ,A)}.
\end{equation}
These mean that $\forall \epsilon  >0$, $\exists T>0$, $\forall t>T$, such that
\begin{equation}
\left \| x(t;\alpha ,A)-x(t;\bar{\alpha} ,A)\right \|<\epsilon  /2,\left \| x(t;\bar{\alpha} ,A)\right \|<\epsilon  /2.
\end{equation}
We have
\begin{equation}
\left \| x(t;\alpha ,A) \right \|< \left \| x(t;\alpha ,A)-x(t;\bar{\alpha} ,A)\right \|+\left \| x(t;\bar{\alpha} ,A)\right \|<\epsilon.
\end{equation}
This is a contradiction.
\par 2.2 Suppose $\exists (\alpha ,A)\in R{(\alpha ,A)}$, $\left | arg\left ( \lambda_{(\alpha ,A)} \right ) \right |=r>\pi /2$. Then there exists a vector $\bar{\alpha}$ whose components are rational numbers, such that
\begin{equation}
r<\left | arg\left ( \lambda_{(\bar{\alpha} ,A)} \right ) \right |<\pi/2,(\bar{\alpha},A)\in S{(\alpha ,A)}\bigcap \delta {(\alpha ,A)}.
\end{equation}
where $\delta{(\alpha ,A)}$ is a neighborhood of $(\alpha,A)$. Eq.(17) means $(\bar{\alpha},A)\notin S{(\alpha ,A)}$. This is a contradiction.
\par Then the conclusion holds.
\endproof
Given $\alpha=(\alpha_1,\alpha_2,\cdots,\alpha_n)$, the corresponding system is of order-certain. Let $z$ be a root of the order-certain system. We call $s$ a critical root and the corresponding $\alpha$ a critical order of system (1) if $s$ satisfies $\left | arg\left ( s \right ) \right |=\frac{\pi }{2}$. Let $s=r\cdot e^{i\frac{\pi}{2}}$, Eq.(3) can be transformed into
\begin{equation}
det\left ( diag\left ( {r^{\alpha_1}e^{i\frac{\pi \alpha_{1} }{2}}},{r^{\alpha_2}e^{i\frac{\pi \alpha_{2}}{2}}},\cdots ,{r^{\alpha_n}e^{i\frac{\pi \alpha_{n}}{2}}} \right )-A \right )=0,
\end{equation}
which can be rewritten as
\begin{equation}
f_1(r,\alpha, A)+i\cdot f_2(r,\alpha, A)=0,
\end{equation}
where $f_1(r,\alpha, A)$ and $f_2(r,\alpha, A)$ are the real part and the imaginary part of Eq.(18), respectively.
\par
The boundary of $S(\alpha,A)$ is contained in the hyper-surface in the parameter space defined by the following parametric equation:
\begin{equation}
\left\{\begin{matrix}
f_1(r,\alpha,A )=0\\
f_2(r,\alpha,A )=0
\end{matrix}\right.,
\end{equation}
where $r$ is the intermediate parameter.
\par
For system (1), the imaginary axis is the critical line which separates stable parameter region and unstable parameter region in complex plane. The way to test the stability of system (1) is to check the distribution of all roots of Eq.(3). When these order parameters are uncertain, the way is hard to effective because the roots may be of parameter-dependence. Based on the above analysis, we transform the critical line into the boundary of $S(\alpha,A)$, so the relationship between parameters and stability can be studied directly.
\par
For system (1), given one point $\alpha=(\alpha_1,\alpha_2,\cdots,\alpha_n)$, the root of the corresponding certain system denoted by $\lambda_{(\alpha,A)}$, which is a continuous function of $\alpha$ and $A$. Denote $\Gamma=(\alpha,A)$, then $\lambda_{(\alpha,A)}$ is shorted as $\lambda_{\Gamma}$.
\begin{theorem}
For system (1), the hyper-surface defined by Eq.(20) decomposes the parameter space into the finite number of connected regions. For any $\Gamma_1,\Gamma_2$ in the same region $\Omega$, we have $Sign(Re(\lambda_{\Gamma_1}))
=Sign(Re(\lambda_{\Gamma_2}))$, where
\begin{equation}
Sign\left (z\right )=\left\{\begin{matrix}
1,z>0\\
0,z=0\\
-1,z<0
\end{matrix}\right..
\end{equation}
\end{theorem}
\proof
Given $\Gamma_1$ and $\Gamma_2$ in the same region $\Omega$, we assume
\begin{equation}
Sign\left (Re\left (\lambda_{\Gamma_1}  \right )  \right )\cdot  Sign\left (Re\left (\lambda_{\Gamma_2}  \right )  \right )=-1.
\end{equation} Taken a curve segment contained in $\Omega$ and connecting $\Gamma_1$ and $\Gamma_2$, because $\lambda_{\Gamma}$ is a continuous function of $\alpha$ and $A$, there exists one point $\Gamma {}'$ on the curve segment $\Gamma_1\Gamma_2$, so that $Sign\left (Re\left (\lambda_{\Gamma{}'}  \right )  \right )=0$. This means there exists one critical root in the region, it is contradicted to that $\Gamma_1$ and $\Gamma_2$ in same region because all critical roots are on the critical line. So the assumption is invalid. The proof is thus completed.
\endproof
\par
Based on Theorem 2, we know that all systems whose parameters belong to the same region have the same stability.  Each region only needs checking one point to determine the stability of the region. This is an effective method for system (1) with uncertain parameters.
\section{Illustrative Examples}
\begin{example} Consider the following characteristic equation which corresponds to a fractional-order system with incommensurate orders:
\begin{equation}
s^\alpha-as^\beta+b=0,
\end{equation}
where $\alpha$ and $\beta$ are different orders, $a$ and $b$ are coefficients.
\end{example}
\begin{figure}[h]
\begin{center}
\includegraphics[width=2.5in]{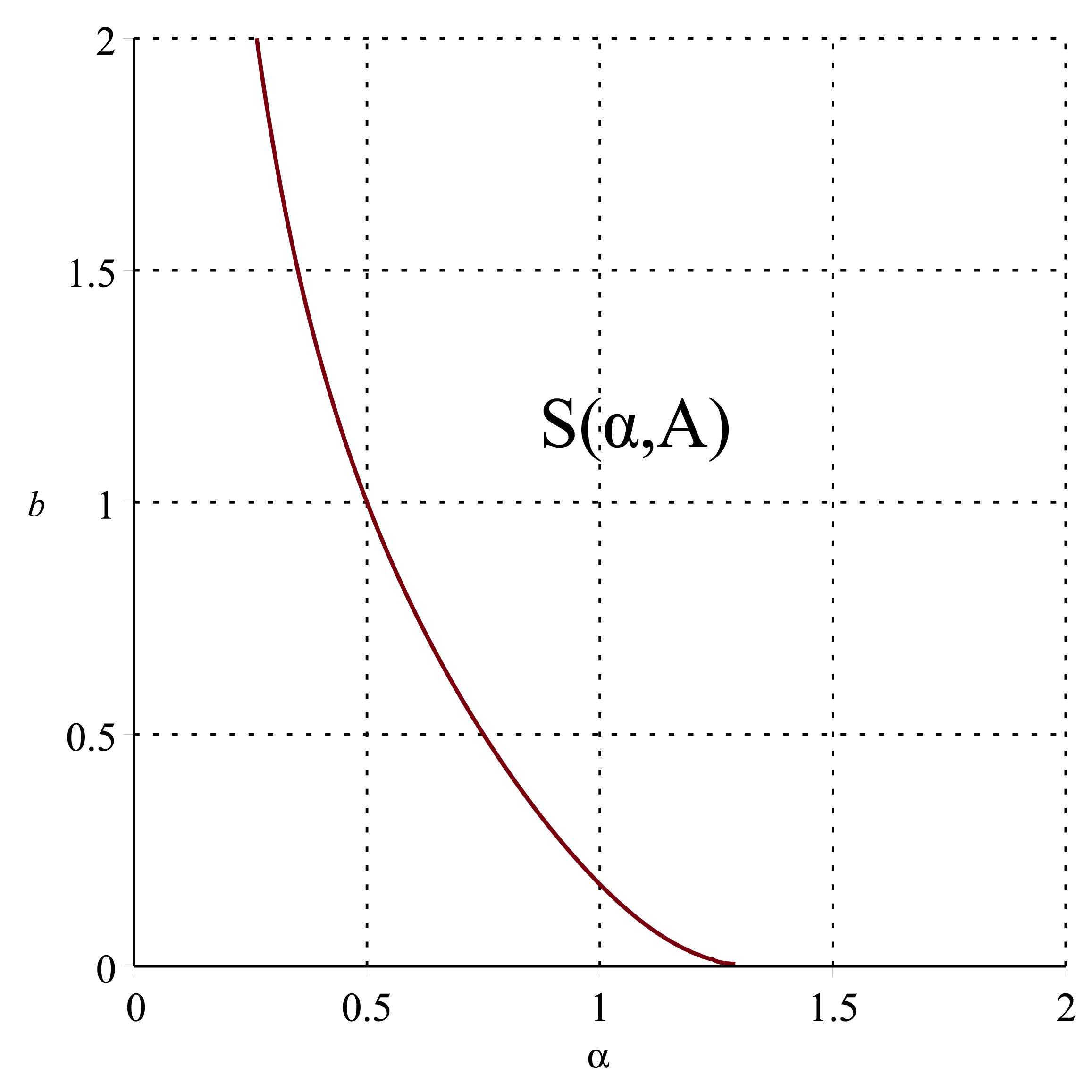}
\end{center}
\caption{The stable parameter region and its boundary of system (23) with $a=-2,\beta=1.5$}
\label{fig1}
\end{figure}
Let $s=r\cdot e^{i\cdot \frac{\pi }{2}}$,
\begin{equation}
\begin{aligned}
s^{\alpha }-as^{\beta }+b &= \left (r\cdot e^{i\cdot \frac{\pi }{2}}  \right )^{\alpha }-a\left (r\cdot e^{i\cdot \frac{\pi }{2}}  \right )^{\beta }+b\\
 &= r^{\alpha }e^{\frac{i\alpha \pi }{2}}-ar^{\beta }e^{\frac{i\beta\pi }{2}}+b\\
 &=r^{\alpha }\left ( cos\left ( \frac{\alpha \pi }{2} \right )+i\cdot sin\left ( \frac{\alpha \pi }{2} \right )  \right )-ar^{\beta }\left ( cos\left ( \frac{\beta \pi }{2} \right )+i\cdot sin\left ( \frac{\beta \pi }{2} \right )  \right )+b\\
 &= r^{\alpha }cos\left ( \frac{\alpha \pi }{2}\right )-ar^{\beta }cos\left ( \frac{\beta \pi }{2} \right )+b+i\cdot \left ( r^{\alpha }sin\left ( \frac{\alpha \pi }{2} \right )-ar^{\beta }sin\left ( \frac{\beta \pi }{2} \right ) \right )\\
 &=0,
\end{aligned}
\end{equation}
which is equivalent to the following parametric equation,
\begin{equation}
\left\{\begin{matrix}
f_1(r)=r^{\alpha }cos\left ( \frac{\alpha \pi }{2}\right )-ar^{\beta }cos\left ( \frac{\beta \pi }{2} \right )+b=0\\
f_2(r)=r^{\alpha }sin\left ( \frac{\alpha \pi }{2} \right )-ar^{\beta }sin\left ( \frac{\beta \pi }{2} \right )=0
\end{matrix}\right.,
\end{equation}
$r^{\alpha}$ and $r^{\beta}$ can be obtained from Eq.(25),
\begin{equation}
r^{\alpha }=\frac{bsin(\alpha \pi/2)}{sin\left (\left (\beta -\alpha   \right )\pi /2  \right )}, r^{\beta } =\frac{-bsin(\alpha \pi/2)}{asin\left (\left (\beta -\alpha   \right )\pi /2  \right )}.
\end{equation}
For $\left (r^{\alpha }  \right )^{\beta }=\left (r^{\beta }  \right )^{\alpha }$, we have
\begin{equation}
\left (\frac{bsin(\alpha \pi/2)}{sin\left (\left (\beta -\alpha   \right )\pi /2  \right )}\right )^{\beta }=\left (\frac{-bsin(\alpha \pi/2)}{asin\left (\left (\beta -\alpha   \right )\pi /2  \right )} \right )^{\alpha }.
\end{equation}
Eq.(27) is the equation of the boundary of stable parameter region, and
\begin{equation}
\begin{aligned}
f(\alpha,\beta) &= \left (\frac{bsin(\alpha \pi/2)}{sin\left (\left (\beta -\alpha   \right )\pi /2  \right )}\right )^{\beta }-\left (\frac{-bsin(\alpha \pi/2)}{asin\left (\left (\beta -\alpha   \right )\pi /2  \right )} \right )^{\alpha }\\
 &= \frac{a^{\alpha }b^{\beta }sin^{\beta }(\alpha \pi/2)sin^{\alpha }\left (\left (\beta -\alpha   \right )\pi /2  \right )-(-b)^{\alpha }sin^{\alpha }(\alpha \pi/2)sin^{\beta }\left (\left (\beta -\alpha   \right )\pi /2  \right )}{a^\alpha sin^{\alpha+\beta }\left (\left (\beta -\alpha   \right )\pi /2  \right )}\\
 &=0
\end{aligned}
\end{equation}
where $\alpha \neq \beta$.
\par Given $a=-2$, $\beta=1.5$, the stable parameter region and its boundary are shown in Fig.1. Our results also can be used to analyze the relationship between uncertain orders and uncertain other parameters.
\begin{remark}
Paper \cite{vcermak2015asymptotic} showed a similar system under the condition $\alpha=2\beta$. As a comparison, the system in Example 1 is of incommensurate orders. Example 1 gives the equation of the boundary of stable parameter region, which includes the relationship among multiple parameters.
\end{remark}
\begin{example} Consider a fractional-order boost converter system as follows:
\begin{equation}
\begin{bmatrix}
\frac{d^{\alpha }I_L}{dt^{\alpha }}\\
\frac{d^{\beta }V_0}{dt^{\beta }}
\end{bmatrix}=\begin{bmatrix}
0 & \frac{-1000}{3}\\
10000 & \frac{-1000}{3}
\end{bmatrix}\begin{bmatrix}
I_L\\
V_0
\end{bmatrix}+\begin{bmatrix}
4000\\
0
\end{bmatrix},
\end{equation}
where $\alpha ,\beta \in (0,2)$, $\alpha \neq \beta$, $I_L$ is inductive current, $V_0$ is output voltage.
\end{example}
\par The characteristic equation of system (29) is
\begin{equation}
s^{\alpha+\beta}+\frac{1000}{3}s^{\alpha}+\frac{10000000}{3}=0,
\end{equation}
the equation of the boundary of stable parameter region is
\begin{equation}
\left (-\frac{1000sin\frac{\alpha \pi }{2}}{3sin\frac{(\alpha +\beta )\pi }{2}}  \right )^{\alpha }=\left (\frac{10000sin\frac{(\alpha +\beta )\pi }{2}}{sin\frac{\beta \pi }{2}}  \right )^{\beta }.
\end{equation}
The stable parameter region and its boundary of system (29) is shown in Fig.2. The boundary decomposes the parameter space into two regions. Based on Theorem 2, we only need to take two points to check the stability of system (29). For point $(\alpha,\beta)=(1,1)$, the corresponding order-certain system is stable, so the corresponding region is the stable region, the other is unstable region.
\begin{figure}[h]
\begin{center}
\includegraphics[width=2.5in]{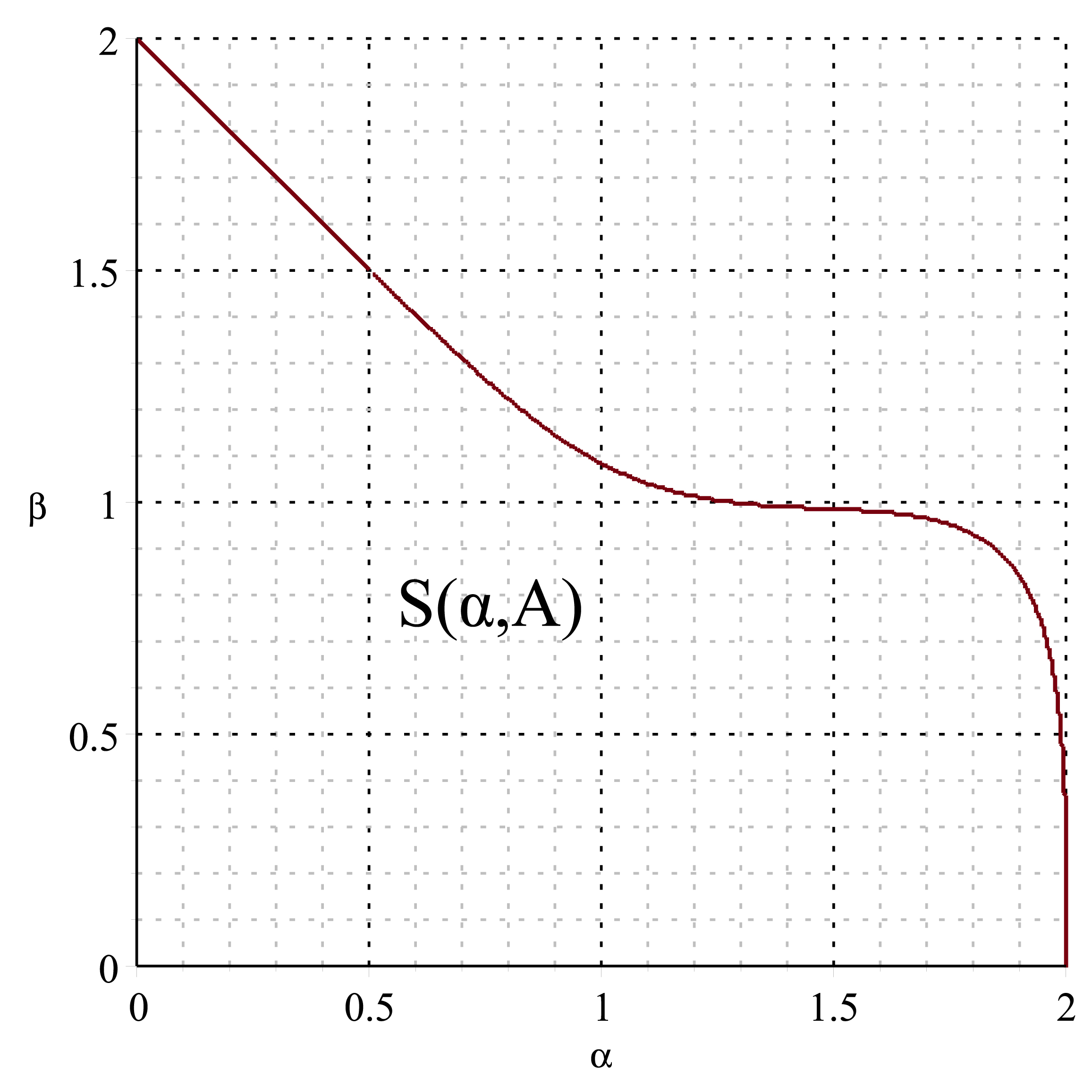}
\end{center}
\caption{The stable parameter region and its boundary of system (29)}
\label{fig2}
\end{figure}
\begin{remark} Paper \cite{brandibur2021exact} discussed a similar system as follows:
\begin{equation}
\left\{\begin{matrix}
\textrm{D}^{q_1}x(t)=0.00001x(t)+y(t)\\
\textrm{D}^{q_2}y(t)=-0.0022x(t)+0.1y(t)
\end{matrix}\right.,
\end{equation}
where $q_1,q_2\in (0,1]$. Based on the results of two-dimensional fractional-order systems, researchers described the stable region based on numerical calculation method. Based on our method, the stable parameter region and its boundary are shown in Fig.3. Taking one point in each region can check the stability easily. Fig.3 is similar to the results in paper \cite{brandibur2021exact}. Our results determine the boundary for the n-dimensional fractional-order system explicitly, then take finite number of points to test which regions are the stable regions. The results are complete which include the results of two-dimensional fractional-order systems.
\begin{figure}[h]
\begin{center}
\includegraphics[width=2.5in]{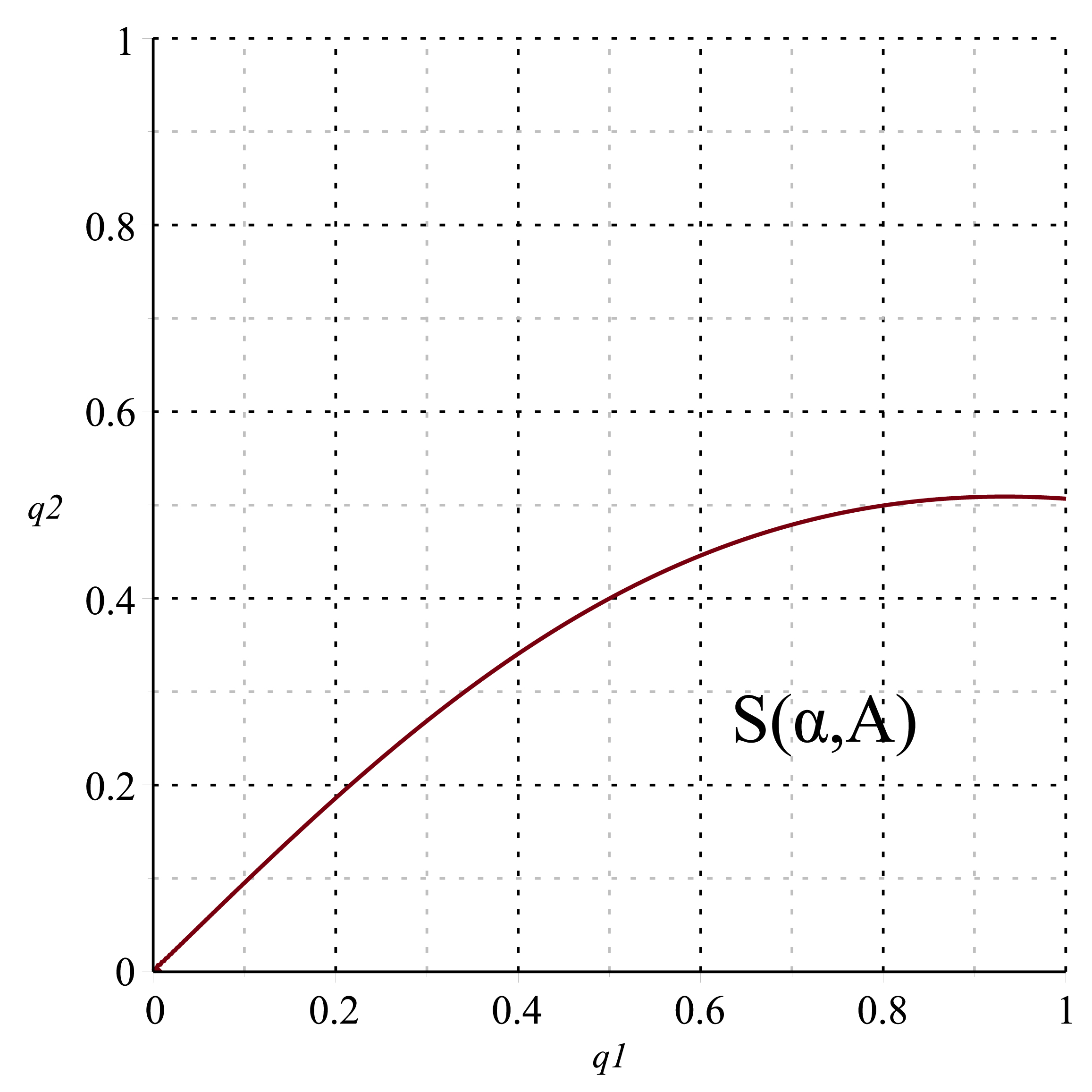}
\end{center}
\caption{The stable parameter region and its boundary of system (32)}
\label{fig3}
\end{figure}
\end{remark}
\begin{remark}
Consider the following fractional-order system with three different orders:
\begin{equation}
\left\{\begin{matrix}
\textrm{D}^{q_1}x_1(t)=x_1(t)-x_2(t)\\
\textrm{D}^{q_2}x_2(t)=2x_2(t)-x_3(t)\\
\textrm{D}^{q_3}x_3(t)=10x_1(t)+10x_2(t)+20x_3(t)
\end{matrix}\right.,
\end{equation}
where $q_{1},q_2,q_3\in\left ( 0,2 \right )$.
\par
The characteristic equation of system (33) is
\begin{equation}
s^{q_1+q_2+q_3}-20s^{q_1+q_2}-2s^{q_1+q_3}-s^{q_2+q_3}+50s^{q_1}+20s^{q_2}+2s^{q_3}-60=0.
\end{equation}
Let $q_3=2q_1$, by using our method, the stable parameter region and its boundary are shown in Fig.4.
\begin{figure}[h]
\begin{center}
\includegraphics[width=2.5in]{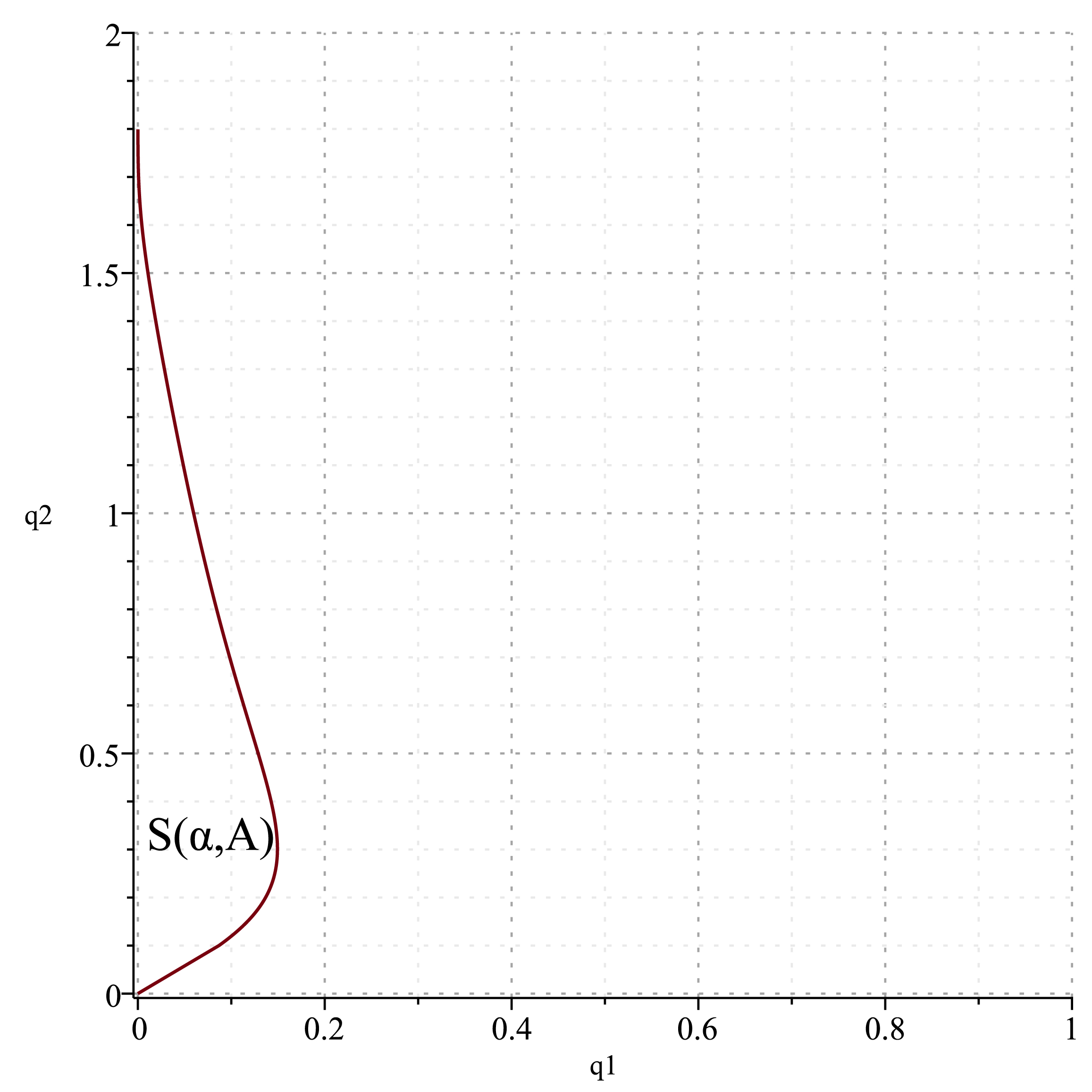}
\end{center}
\caption{The stable parameter region and its boundary of system (33)}
\label{fig4}
\end{figure}
\end{remark}
\begin{example} Consider the characteristic equation:
\begin{equation}
s^{\alpha_1+\alpha_2}+12s^{\alpha_1}+34=0.
\end{equation}
\end{example}
When the LCM of the denominators of $\alpha_1$ and $\alpha_2$ is big, for example, $\alpha_1=993/1000,\alpha_2=997/1000$, then Eq.(35) can be transformed to an integer-order characteristic equation:
\begin{equation}
s^{1990}+12s^{993}+34=0,
\end{equation}
the number of roots is too large to calculate. Our method is suitable for this situation. The equation of the boundary of stable parameter region is
\begin{equation}
\left (-\frac{12sin\frac{\alpha \pi }{2}}{sin\frac{(\alpha +\beta )\pi }{2}}  \right )^{\alpha }=\left (-\frac{17sin\frac{(\alpha +\beta )\pi }{2}}{6sin\frac{\beta \pi }{2}}  \right )^{\beta }.
\end{equation}
The stable parameter region and its boundary are shown in Fig.5. We see that the boundary decomposes the parameter region into three regions. For point $(\alpha,\beta)=(1,1)$, the corresponding system is stable, so the corresponding region is the stable parameter region. And the others are unstable parameter regions, which have different distributions of roots.
\begin{figure}[h]
\begin{center}
\includegraphics[width=2.5in]{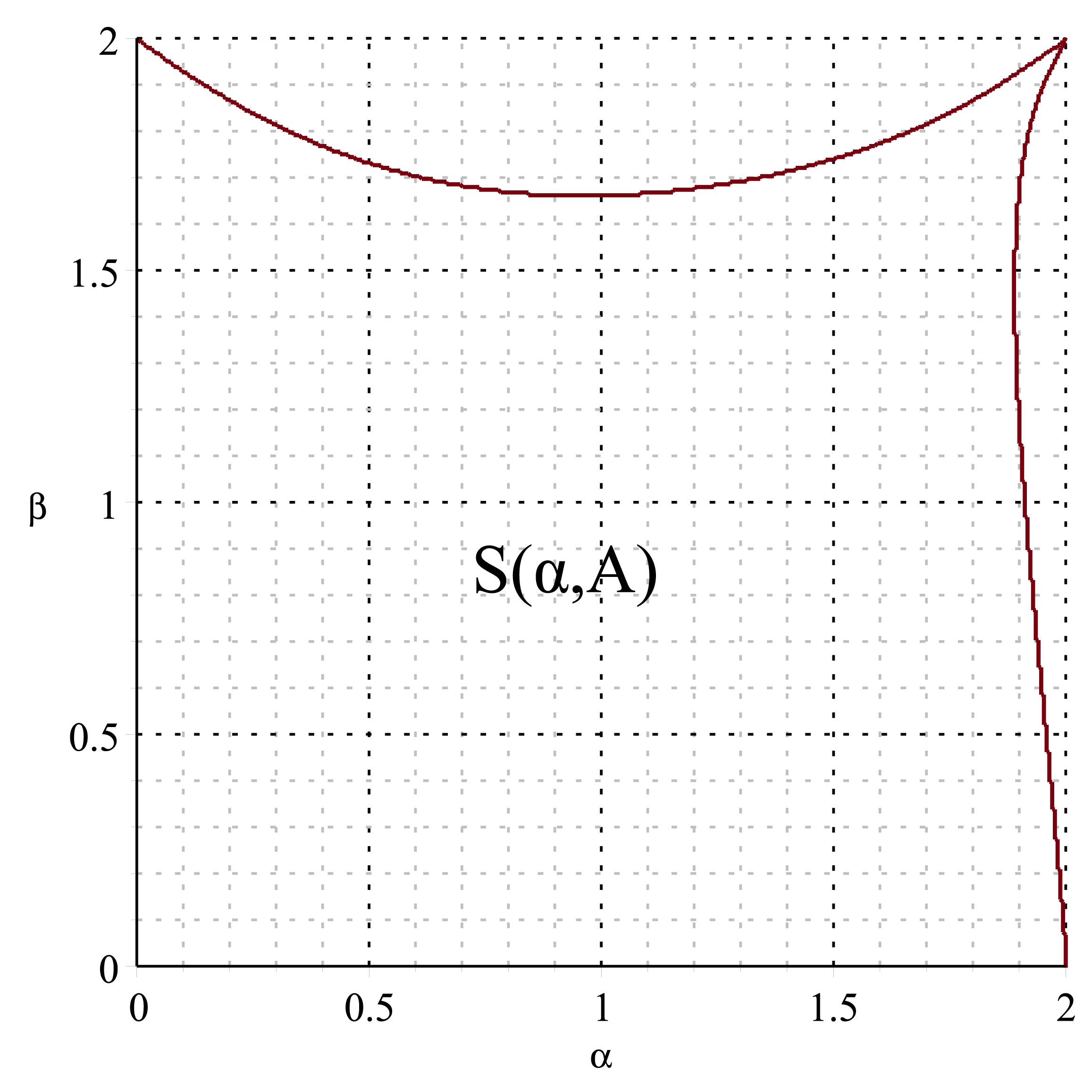}
\end{center}
\caption{The stable parameter region and its boundary of system (35)}
\label{fig5}
\end{figure}
\begin{remark}
The existing researches analyzed fractional-order systems with incommensurate orders by using Lemma 1. The number of roots is large in some examples. For instance,
paper \cite{daftardar2010chaos} showed the following characteristic equation of fractional-order Liu system:
\begin{equation}
\lambda^{129}+\lambda^{100}+2.5\lambda^{89}+7.5\lambda^{50}-1.7764\times 10^{-15}\lambda^{39}+50=0.
\end{equation}
This equation has 129 roots. That system has two different orders. By using our method, the computation is less and more information between orders can be obtained.
\end{remark}
\begin{remark}
Consider the existing systems in papers \cite{liao2011robust,li2012robust}. If orders are interval value, then Lemma 1 will fail.
\end{remark}
\section{Conclusions}
Fractional-order systems with incommensurate orders are considered in this paper. In most existing papers, researchers considered order parameters are rational numbers which limited discussion of some cases just like order parameters are uncertain parameters or interval parameters. In this paper, we discuss the problem of stability on n-dimensional systems which have real number order parameters. The stability criterion is determined by strict proof. Based on the results, we determine the boundary of stable parameter region, which decomposes the parameter space into finite number of connected regions. All systems whose parameters belong to the same region have the same stability. Each region only needs checking one point to determine the stability of the region. The results are complete and include the relationship among different orders. Several examples demonstrate the effectiveness of our method.
\section*{Acknowledgment}
This work was partially supported by the National Natural Science Foundation of China
under Grants No. 12171073.
\section*{References}
\bibliographystyle{elsarticle-num-names}
\biboptions{square,numbers,sort&compress}
\bibliography{elsarticle-template}

\end{document}